\documentclass[11pt]{amsart}
\usepackage{amsmath}
\usepackage[all, cmtip]{xy}
\ifx\pdfoutput\undefined 
\usepackage{graphicx}
\else
\usepackage[pdftex]{graphicx}
\usepackage{epstopdf}
\fi
\usepackage{verbatim}

\newtheorem{theorem}{Theorem}[section]
\newtheorem{proposition}[theorem]{Proposition}
\newtheorem{lemma}[theorem]{Lemma}

\theoremstyle{definition}
\newtheorem{definition}[theorem]{Definition}

\theoremstyle{remark}
\newtheorem{remark}[theorem]{Remark}
 
\numberwithin{equation}{section}

\newcommand{\PB}[2]{\{#1,#2\}}

\newcommand{\h}{\hbar}

\newcommand{\bbC}{{\mathbb C}}

\newcommand{\bbZ}{{\mathbb Z}}
\newcommand{\bbR}{{\mathbb R}}

\newcommand{\calR}{{\mathcal R}}

\newcommand{\zbar}{\overline{z}}

\DeclareMathOperator{\tr}{Tr}

\newcommand{\av}{\mbox{\tiny{ave}}}

\begin{document}

\openup1pt

\title{Band invariants for perturbations of the harmonic oscillator}
\author{V. Guillemin}\address{Department of Mathematics\\
Massachusetts Institute of Technology \\Cambridge, MA 02139 \\
USA}\thanks{V. Guillemin is supported in part by NSF grant
DMS-1005696.}\email{vwg@math.mit.edu}
\author{A. Uribe}\address{Department of Mathematics \\
University of Michigan \\ Ann Arbor \\ MI \\ 48109 \\ USA}
\thanks{A. Uribe is supported in part by NSF grant DMS-0805878.}
\email{uribe@umich.edu}
\author{Z. Wang}
\address{Department of Mathematics \\
University of Michigan \\ Ann Arbor \\ MI \\ 48109 \\ USA} \email{wangzuoq@umich.edu}

\date{}

\begin{abstract}

We study the direct and inverse spectral problems for semiclassical operators of the form
$S = S_0 +\h^2V$, where $S_0 = \frac 12 \Bigl(-\h^2\Delta_{\bbR^n} + |x|^2\Bigr)$ is the 
harmonic oscillator and $V:\bbR^n\to\bbR$ is a tempered smooth function.
We show that the spectrum of $S$ forms eigenvalue clusters as $\h$ tends to zero, and
compute the first two associated ``band invariants".  We derive several inverse spectral results
for $V$, under various assumptions.  In particular we prove that, in two dimensions,
generic analytic potentials that are even with respect to each variable are spectrally determined
(up to a rotation).
\end{abstract}

\maketitle
\tableofcontents

\section{Introduction}

Consider the semiclassical harmonic oscillator  
\[ 
  S_0 = \frac{1}{2}\Bigl(-\h^2 \Delta_{\mathbb R^n} + |x|^2 -\h\,nI\Bigr) 
\]
acting on $L^2(\mathbb R^n)$, where, to simplify the notation throughout the paper, 
we have subtracted the ground state energy.
The spectrum of $S_0$ consists of the eigenvalues
  \begin{equation}\label{1.1}
E_j(\h) = \h j , \quad j=0, 1, 2, \cdots,
  \end{equation}
with multiplicity
\[
m_j = {n+j-1\choose n-1}.
\]

We will study perturbations of $S_0$ of the form
  \begin{equation}\label{GenPer}
  S = S_0 + \h^2 A,
  \end{equation}
where $A$ is a self-adjoint semiclassical
pseudodifferential operator of order zero. Since the distance between consecutive eigenvalues of 
$S_0$ is $O(\h)$, if $A$ is $L^2$ bounded for $\h$ sufficiently small the spectrum of $S$ 
forms non-overlapping clusters of eigenvalues 
  \[ 
 \{E_{j,k};\;  k=1, \cdots, m_{j}\}, \quad j=0, 1, \cdots, 
  \]
defined by the condition 
  \[ 
  |E_{j,k}(\h)  -E_j(\h) | = O(\h^2). 
  \]
As we will see the $L^2$ boundedness assumption of $A$ can be dropped, provided one looks 
at the
spectrum of $S$ locally:  For all $A$ in standard symbol classes the spectrum of $S$ restricted 
to compact intervals forms spectral clusters for small enough $\h$.

By analyzing these spectral clusters using wave-trace
techniques, one generates a collection of spectral invariants which are
known as ``band invariants" and have been studied in detail in various
contexts in \cite{CV}, \cite{Wei}, \cite{Gu1}, \cite{Gu2}, \cite{Doz}, \cite{Ur}, \cite{UrV}, 
\cite{Zel}, and elsewhere.
The main topic of this paper will be computations of these invariants for
the harmonic oscillator. These computations will be valid in arbitrary
dimension, but in dimension two we will be able to extract from them a
number of new inverse results. In particular, for the Schr\"odinger-like
operator $S_0 + \h^2 V$, we will show that if V is real analytic and
symmetric in its  $x_1$ and $x_2$ coordinates, and if, in addition, the two
eigenvalues of the Hessian of $V$ at the origin are distinct, then $V$ is
spectrally determined. In fact we will show more generally that a result of
this nature is true for real analytic semi-classical potentials, $V =
V(x,\h)$, and for $C^\infty$ semi-classical potentials  provided that $V(x,0)$
is a quadratic form with distinct eigenvalues. (See section 6 for
details.)

Perturbations of the form $S_0 + \h^{1+\delta} A$, $\delta >0$ and $A$ a bounded
zeroth order pseudodifferential operator, have been 
studied recently by D. Ojeda-Valencia and C. Villegas-Blas , \cite{O-VV-B}.  Working in 
Bargmann space they derive the first band invariant for such perturbations. 

\section{The Clustering Phenomenon}

We will be working with semiclassical pseudodifferential operators, with amplitudes
in the symbol classes ${\mathcal O}(\langle x,\xi\rangle^{m})$, and 
using the Weyl calculus (see \cite{Mar} for definitions).

\subsection{Existence of Spectral Clusters}
Let $S$ be as in (\ref{GenPer}). 
The following theorem is a special case of a result by B. Helffer and D. Robert
in \cite{HR} (Theorem 3.9).  
It states that, if we only look at eigenvalues contained in a fixed ``window", 
then as $\hbar \to 0$, the eigenvalues of $S$ form clusters around the
unperturbed eigenvalues. 
We include a proof here for completeness, and to set up the
computation of the band invariants.

\begin{theorem}\label{Clusterexistence}
For any compact interval $I \subset \bbR$ there exists  $\hbar_0>0$ small enough such 
that for $0<\hbar<\hbar_0$, the spectrum of $S$  in $I$ consists of eigenvalues that cluster 
near the eigenvalues of $S_{0}$ within a distance $O(\hbar^2)$.  More precisely, given
$I$
\[
\exists C>0,\ \h_0>0 \quad\text{such that}\quad \forall h\in (0,\h_0]
\]
\begin{equation}\label{Cluster}
\text{Spec }(S)\cap I \subset \bigcup_{j\in\bbZ_+} \left[\h j - C\h^2\,,\, \h j + C\h^2\right]  .
\end{equation}
\end{theorem}
\begin{proof}
Define the time-dependent operator $R(t)$ by the identity
\begin{equation}\label{defR}
\; e^{i\h^{-1}t S_{0}}e^{-i\h^{-1}t S} = I +\h R(t).
\end{equation}
The left hand side of this identity is a semiclassical pseudodifferential operator (a composition 
of semiclassical FIOs with
inverse canonical relations), so $R(t)$ is a semiclassical pseudodifferential operator. In what 
follows we let 
\[R=R(2\pi) = \frac 1{\h}(e^{-2\pi i \h^{-1}S} - I).\]

Now let $\chi \in C_0^\infty(\mathbb R)$ be a cut-off function which is identically equal to one 
in an open neighborhood of $I$, and define 
\[R_\chi = \chi(S)R.\]
This is a semiclassical pseudodifferential operator, and in the expansion of its symbol
\[ r^\chi(x, p, \h)  \sim \sum h^j r^\chi_j(x, p),\] 
one can see that each $r_j^\chi$ has compact support contained in 
$H_0^{-1}(\mathrm{supp}(\chi))$, where 
\begin{equation}\label{Hzero}
H_0(x, p) = \frac 12 (|x|^2+|p|^2)
\end{equation}
is the principal symbol of $S_0$. 
It follows from the Calder\'on-Vaillancourt theorem that $R_\chi$ a bounded 
operator in $L^2$, and the bound is \emph{uniform} in $\h$. 
Thus, for $\h$ sufficiently small the spectrum of $I+\h R_\chi$ is contained in a neighborhood
of 1 away from
the origin in $\bbC$, and we can define
\begin{equation}\label{defW}
W = \frac{i}{2\pi\h} \, \log\left(I + \h R_\chi\right),
\end{equation}
which is a bounded
pseudodifferential operator of order zero.  Define now
\[
S_\chi := S -\h^{2} W.
\]
Note that all the operators in this paragraph are functions of $S$, and therefore commute
with each other.

It follows from the definitions that
\[
e^{-2\pi i \h^{-1} S_\chi} = (I+\h R)\, (I + \h R_\chi)^{-1}.
\]
Now let $\Lambda\in I$ be an eigenvalue of $S$ with normalized eigenvector $\psi$.  Since 
\[
\chi(S)(\psi) = \psi,
\]
it is clear that
\[
e^{-2\pi i \h^{-1} S_\chi}(\psi) = \psi,
\]
which is to say that $S_{\chi}\psi = \hbar j\,\psi$ for some integer $j$.
From this and the fact that $S = S_{\chi}+\h^{2} W$ it follows that
\begin{equation}\label{bandseq}
\Lambda = \h j + \h^{2}\mu,\quad\text{where}\quad W\psi = \mu\psi.
\end{equation}
Therefore 
\[
|\Lambda-\h j| \leq C\h^{2},
\]
where $C$ is a uniform $L^{2}$ bound on the norm of $W$ for $\h$ small enough.  
This completes the proof.
\end{proof}

Note that the intervals on the right-hand side of (\ref{Cluster}) are pairwise
disjoint for $\h$ small enough, since their centers are a distance $\h$ apart and their lengths
are $O(\h^{2})$.  The eigenvalues of $S$ in any one of those disjoint intervals must be of the 
form
\begin{equation}\label{specShifts}
E_{j,k} = E_j + \h^{2}\mu_{j,k}\,
\end{equation}
where the $\mu_{j,k}$ are the eigenvalues of $W$ 
restricted to the eigenspace of $S_\chi$ corresponding to $E_j =\hbar j$.
Note that the $\mu_{j,k}$ are uniformly bounded.

\begin{remark}
The previous theorem holds for more general \emph{Zoll operators}.  
Namely, let $Q$ be a semiclassical pseudodifferential operator such that Spec$(Q) \subset \hbar 
\mathbb Z$.  (This implies that the Hamilton flow of its principal symbol is $2\pi$ periodic.)
Let $A$ be a zeroth order semiclassical
pseudodifferential operator whose symbol lies in ${\mathcal O}(\langle x,\xi\rangle^{m})$, and 
consider the perturbation of $Q$
\[
P = Q + \hbar^2 A.
\]
Then for any compact interval $I \subset \mathbb R$ the eigenvalues of $P$ in $I$ cluster near 
$\hbar \mathbb Z$, for all $\hbar$ small enough.  
\end{remark}

\medskip
\noindent
{\em Clusters for potentials of quadratic growth.}

In the case that the perturbation $A$ is a ``multiplication by a potential function $V(x)$" 
operator, there is a much simpler and direct proof of the previous result provided $V$ is of no 
more than quadratic growth at infinity. 
Moreover, in this case one can see the spectral clustering appearing uniformly
in all of $\bbR$, instead only locally uniformly in compact intervals. 
In other words, there exists $C>0$ and $\h_0>0$ such that for all $\h<\h_0$, 
\begin{equation}\label{ClusterQ}
\text{Spec}(S_0+\h^2 V) \subset \bigcup_{k \in \bbZ_+} \left[ \h k - C\h^2, \h k + C \h^2 
\right].
\end{equation}
Here is a sketch of the proof: If $V$ is of no more than quadratic growth at infinity, 
then one can find constants $C_1>0, C_2>0$ such that 
\[ -C_1 |x|^2 - C_2 <V(x)  < C_1 |x|^2 +C_2.\] 
On the other hand, by directly computations, one can show that the eigenvalues of 
\[S_0+\h^2(\pm C_1 |x|^2 \pm C_2)  =
- \frac 12  \hbar^2\Delta + \frac 12 |x|^2(1\pm 2\h^2 C_1) \pm \h^2 C_2\]
are precisely $\h \sqrt{1\pm 2\h^2 C_1} (j+\frac n2) \pm  \h^2 C_2$, which are within 
$O(\h^2)$ distance of the unperturbed eigenvalues (\ref{1.1}).
Now (\ref{ClusterQ}) follows from a min-max argument.

\begin{remark}
Analogous arguments show that the same result holds for a ``semiclassical potential'',
\[
A=V_0(x)+\hbar V_1(x)+ \hbar^2 V_2(x) + \cdots,
\] 
where all the $V_j$ are of functions with uniform quadratic growth :  
$|V_j(x)|<C_1+C_2 |x|^2$ (for some constants $C_1, C_2$ independent of $j$).
\end{remark}

\subsection{Cluster projectors}
Let $I\subset \bbR$ be a compact interval.  Keeping the notations of 
Theorem \ref{Clusterexistence}, let $E\in\text{Int} (I)$, $E>0$,
and restrict the values of $\h$ to the sequence
\begin{equation}\label{hsequence}
\h = \frac{E}{N},\quad N=1,\ 2,\ldots .
\end{equation}
For each such $\h$ small enough,
\begin{equation}\label{indcluster}
\text{Spec}(S)\cap [E-C\h^2\,,\,E+C\h^2] = \{ E+h^{2 }\mu_{N,k}\;;\; k=1,\ldots , m_{N}\}
\end{equation}
is a single cluster (see (\ref{specShifts})).

\begin{theorem} \label{ClustProj}
For each $\h$ as in (\ref{hsequence}), 
let $\Pi_{E}^{N}$ be the orthogonal projection from $L^{2}(\bbR^{n})$ onto the span
of the eigenfunctions of $S$ corresponding to eigenvalues in the individual
cluster (\ref{indcluster}).  Then the family $\{\Pi_{E}^{N}\}_{N}$ is a semiclassical Fourier 
integral operator associated with the canonical relation
\[
\{ (\bar x, \bar y) \in \mathbb R^{2n} \times \mathbb R^{2n}\ |\ H_0(\bar x) = E = H_0(\bar y), 
\bar x \mbox{\ lies on the same\ }H_0 \mbox{\ orbit as\ } \bar y\}.
\]
\end{theorem}
\begin{proof}
 Take a closed interval $I_{1}$ containing $E$ in its interior and included in $\text{Int}(I)$,
  and a cutoff function $\chi_1$ such that $\chi_1=1$ on $I_1$ and $\chi_1=0$ 
  outside $I$.   Define
 \[
 U(t) = \chi_1(S_\chi) e^{-2\pi i\hbar^{-1}tS_\chi}.
 \]
Then one can check that (recalling that $\hbar = E/N$),
\begin{equation}\label{fioProof}
\frac 1{2\pi} \int_0^{2\pi} U(t) e^{it\hbar^{-1}E}dt 
 =  \chi_1(S_\chi) \Pi_E^{N}.
\end{equation}
Now $U(t)$ is a Fourier integral operator associated with the graph of the principal
symbol of $S_{\chi}$.  The conclusion of the theorem follows from this and the 
calculus of FIOs  applied to the left-hand side of (\ref{fioProof}).
\end{proof}

\section{Band invariants}

\subsection{The first band invariant} $ $

According to equation (\ref{defR}), the time derivative of $R(t)$ is given by
\begin{equation}\label{dotR}
\dot{R}(t) = \frac{i}{\h} [S_0, R(t)] - iA - i\h AR(t).
\end{equation}
It follows that its principal symbol, $r_0(t)$, satisfies
\[
\dot{r_0}(t) = \PB{H_0}{r_0} - ia \quad\text{and}\quad r_0|_{t=0} = 0,
\]
where  $a$ is the principal symbol of $A$.  The solution to this problem
is
\[
r_0(t, x,p) = -i\int_0^t a(\phi_s(x,p))\, ds,
\]
where $\phi_s$ is the Hamilton flow of $H_0$.
 
In what follows we will denote, for any function $b$ on $\bbR^{2n}$,
\begin{equation}\label{defAve}
b^{\av}(x,p) := \frac 1{2\pi}\int_0^{2\pi} b(\phi_s(x,p))ds.
\end{equation}
With this notation, the principal symbol of 
\[
W = \frac{i}{2\pi\h} \, \log\left(I + \h R_\chi\right)
\]  
is $a^{\av}$.

The first band invariant arises as follows:

\begin{theorem}\label{BandInv1Thm}
For any $f \in C_0^\infty(\mathbb R)$ and $\varphi \in C^\infty(\mathbb R)$, the integral 
  \begin{equation}\label{BandInv1}
  \int f\left( \frac {|x|^2+|p|^2}2\right) \varphi(a^{\av})\, dxdp
  \end{equation}
is a spectral invariant for the family of operators (\ref{GenPer}). 
\end{theorem}
\begin{proof}
Choose the interval $I$ in Theorem \ref{Clusterexistence} 
such that $\mathrm{supp}f \subset\subset I$. Then 
  \[ 
  \sum_j f(E_j) \sum_{k=1}^{m_j} \varphi(\mu_{k,j}) = \tr\left[ f(S_\chi)\varphi(W) \right] . 
  \]
As $\h \to 0$, the trace on the right-hand side has an asymptotic expansion with leading 
coefficient (\ref{BandInv1}).
\end{proof}

As a consequence, one can see that for any $E >0$, the numbers 
\begin{equation}\label{LocBandInv1}
   \int_{|x|^2 +|p|^2= E} \varphi(a^{\av}) d\lambda 
\end{equation}
are spectral invariants of the semiclassical family of operators $S = S_0 + \hbar^2 A$, where $d
\lambda$ is the (normalized) 
standard Lebesgue measure on the sphere $|x|^2+|p|^2=E$.  
In fact, these numbers arise as Szeg\"o limits:
\begin{theorem}
Fix $E>0$ and $\varphi\in C^{\infty}({\mathbb R})$.  Then, as $N\to\infty$
\begin{equation}\label{Sz}
\frac{1}{m_{N}}\sum_{k=1}^{m_{N}}\varphi(\mu_{N,k}) = 
\int_{|x|^2+|p|^2=E} \varphi(a^{\av}) d\lambda + O(1/N).
 \end{equation}
\end{theorem}
The proof is standard:  The left hand side of (\ref{Sz}) is the normalized
trace of $\varphi( W )\Pi^{N}_{E}$.  The asymptotic behavior of
this trace can be computed symbolically given that the cluster projectors
are FIOs (Theorem \ref{ClustProj}).
 
\subsection{Some properties of the averaging procedure} $ $

For future reference we gather here some properties of  the averaging procedure,
restricted to functions $V$ of the position variable $x$ alone.  For such functions,
letting $z=x+ip$, one can rewrite the integral (\ref{defAve}) as  
\[
V^{\av}(z) = \frac 1{2\pi} \int_0^{2\pi} V\left(\frac{e^{it}z+e^{-it}\bar z}2\right) dt.
\]
From this integral formula one can derive many properties of $V^{\av}$ as a function on $
\bbC^n$. Here is a very incomplete list: 
\begin{enumerate}
\item $V^{\av}$ is always an even function. Moreover, one can regard $V^{\av}$ as a function 
defined on $S^{2n-1}/S^1 = \mathbb{CP}^{n-1}$. 

\item  
corresponding complex sphere $|z|^2 = r^2$. As a consequence, 
If $V_1(x) = V_2(x)$ on a ball $|x|^2 \le r^2$, then $V_1^{\av} = V_2^{\av}$ on the complex 
ball $|z|^2 \le r^2$. 

\item If $V(x)$ is an odd function, then $V^{\av} \equiv 0$. (In later applications, we will 
mainly focus on the $V$'s that are even in \emph{each} variable. )

\item If $V(x)$ is of polynomial growth, 
\[|V(x)| \le C_1+C_2|x|^m,\] 
so is $V^{\av}(z)$, i.e. 
\[|V^{\av}(z)| \le C_1 + C_2 |z|^m.\] 
More generally, if $|V(x)| \le h(|x|)$ for some \emph{increasing} function $h$, then $|V^{\av}
(z)| \le h(|z|)$. 
\begin{enumerate}
\item[(4a)] In general, one cannot drop the ``increasing" assumption on $h$ above. 
For example, consider $n=1$ and $V(x) = \frac 1{1+|x|^2}$, then $V(x) = O(\frac 1{|x|^2})$ 
as $x \to \infty$, and
\[\aligned
V^{\av}(z)  &= \frac 1{2\pi} \int_0^{2\pi} \frac 1{1+(\frac {e^{it}z+e^{-it}\bar z}2 )^2} dt \\
&=\frac 1{2\pi} \int_0^{2\pi} \frac 1{1+\frac{|z|^2}2 (1+\cos(2t))} dt.
\endaligned\]
One can check  $V^{\av}(z) \neq O(\frac 1{|z|^2})$.  
\item[(4b)] However, if $V(x) = o(1)$, then $V(z) = o(1)$. We will leave the proof as an 
exercise.
\item[(4c)] Similarly, if $V(x) \ge C |x|^N$ for $|x|$ large enough, then $V^{\av}(z) \ge C'|z|^N
$ for $|z|$ large enough. 
\end{enumerate}
\item If $V(x)$ is homogeneous of degree $m$ in $x$, then $V^{\av}(z)$ is also 
homogeneous of degree $m$ in $z$.  

\item If $V(x) = x^\alpha = x_1^{\alpha_1} \cdots x_n^{\alpha_n}$ is a monomial, then
$V^{\av}$ is zero unless $|\alpha|$ is even, in which case
\begin{equation}\label{AveMon}
\aligned
V^{\av}(z) &= \frac 1{2\pi} \int_0^{2\pi}\left( \frac {e^{it}z_1 + e^{-it}\bar z_1}2\right)^
{\alpha_1} \cdots \left( \frac {e^{it}z_n + e^{-it}\bar z_n}2\right)^{\alpha_n} dt \\
& = \frac 1{2^{|\alpha|}} \sum_{j_1+ \cdots + j_n = \frac{|\alpha|}2 \atop \alpha_r\ge j_r \ge 0} 
{\alpha_1 \choose j_1}\cdots {\alpha_n \choose j_n} z_1^{j_1}\bar z_1^{\alpha_1-j_1} \cdots 
z_n^{j_n}\bar z_n^{\alpha_n-j_n}.
\endaligned
\end{equation}
 \item  $V^{\av}$ is identically zero if and only if $V$ is an odd function.  We will
 prove this in the $n=2$ case in \S 5.  The case $n>2$ reduces to the two dimensional
 case,  as the trajectories of the harmonic oscillator in any 
 dimension lie on two dimensional planes.
\end{enumerate}

\subsection{Higher band invariants for perturbations by potentials} $ $

To derive higher order band invariants, one must look at higher order terms in the symbol of $W
$. Recall that by the semiclassical Weyl calculus, the Weyl symbol of a composition $PQ$ is 
\[\aligned 
a \#_\h b & \sim e^{\frac{i\hbar}2 (D_x D_q - D_y D_p)}(a(x,p)b(y,q))|_{y=x,q=p}\\
& =\sum \hbar^j B_j(a, b)
\endaligned\]
for some bi-differential operators $B_j$, where
\[
B_0(a, b) = ab, \qquad \mathrm{and}\quad B_1(a,b) = \frac{1}{2i}\PB{a}{b}.
\] 
In particular, the Weyl symbol of the commutator $\frac i{\hbar}[S_0, R]$ is given by 
\[\PB{H_0}{ r} + \hbar^2 \PB{H_0}{r}_3 + \hbar^4 \PB{H_0}{ r}_5 + \cdots,\]
where $\{\ , \ \}$ is the Poisson bracket, and $\{\ , \ \}_k$ is the ``higher poisson bracket"
\[\{a, b\}_k = \sum_{|\alpha|+|\beta|=k} \frac 1{\alpha ! \beta !}(-1)^{|\alpha|} \partial_p^\beta 
\partial_x^\alpha a \partial_p^\alpha \partial_x^\beta b. \]

Let us now assume that the operator $A$ is the operator ``multiplication by a potential function 
$V(x)$", i.e. 
\begin{equation}\label{PotPer}
S = S_0 + \hbar^2 V.
 \end{equation}
Since $H_0=\frac 12 (|x|^2+|p|^2)$ is quadratic in both $x$ and $p$,  all higher order Poisson 
brackets of $H_0$ with $r(t)$ vanish. If we write
\[r(t, \hbar, x, p) \sim r_0(t, x, p) + \hbar r_1(t, x, p) + \hbar^2 r_2(t, x, p) + \cdots,\]
and if  abbreviate $r_k(t,x,p)$ to $r_k(t)$, then from (\ref{dotR}) we get 
\[
\sum_k \hbar^k \dot r_k(t) = \sum_k \hbar^k \PB{H_0}{r_k(t)} - i V - i \sum_{j,l} \hbar^{j+l
+1} B_l(V, r_j)
\]
with initial conditions $r_k(0)=0$. In particular, 
\[
\dot r_1(t) = \PB{H_0}{r_1} - i V r_0  \quad\text{and}\quad r_1|_{t=0} = 0.
\]
It is easy to check that the solution to this equation is simply 
\[r_1(t, x, p) = \frac 12 r_0^2(t, x, p).\]
Similarly the equation for $r_2$ is 
\[
\dot r_2(t) = \PB{H_0}{r_2} - i (V r_1+ B_1(V, r_0)) \quad\text{and}\quad r_1|_{t=0} = 0,
\]
and its solution is
\[
r_2(t,x,p) = \frac 16 r_0^3 + 
\frac{i}{2} \int_0^t \int_0^{t-s} \PB{V(\phi_s(x, p))}{V(\phi_{s+u}(x,p))}\ duds.
\]
In general, for any $k \ge 1$, the $r_k(t, x, p)$ is the solution to the problem
\begin{equation}\label{theProblem}
\dot r_k(t) =  \PB{H_0}{r_k} - i \sum_{l=0}^{k-1}B_l(V, r_{k-1-l})  \quad\text{and}\quad r_k|
_{t=0} = 0.
\end{equation}
Note that if we let 
\[
g_k(t) = \sum_{l=0}^{k-1}B_l(V, r_{k-1-l}),
\]
then $g_{k}$ depends only on $r_0, \cdots, r_{k-1}$, in particular it is 
independent of $r_k$. So the problem (\ref{theProblem}) 
can be solved via Duhamel's principle iteratively:
\[
r_k(t, x, p) = -i\int_0^t g_k(t-s, \phi_s(x, p)) ds.
\]
Replacing $t$ by $2\pi$, we get the higher-order 
symbols of $R$, and thus the symbol of $R_\chi$ in the region $H_0^{-1}(I)$. 

To calculate the higher-order symbols of $W$, note that 
\[\aligned 
W & = \frac{i}{2\pi\h} \, \log\left(I + \h R_\chi\right) \\
& =  \frac {i}{2\pi} \left( R_\chi -\frac{\hbar}2 R_\chi^2 +  \frac{\hbar^2}3 R_\chi^3 + \cdots 
\right) .
\endaligned\] 
So if we let 
\[w  
\sim w_0 + \h w_1 + \cdots
\] 
be the full symbol of $W$, then $w_0=V^{\av}$, and
\[w
_k  = \frac 1{2\pi i} \sum_{j=1}^{k+1} \frac{(-1)^{j+1}}j \sum_{l_1+\cdots+l_{j-1}+m_1+
\cdots+m_j=k-j+1} B^{l_1, \cdots, l_{j-1}}(r_{m_1}, \cdots, r_{m_j}),
\]
where $r_j$'s are symbols of $R_\chi$ calculated above, all indices
$l_j, m_j$'s are nonnegative, and
\[ 
B^{l_1, \cdots, l_{j-1}}(r_{m_1}, \cdots, r_{m_j}) = B_{l_{j-1}}(B_{l_{j-2}}(\cdots (B_{l_1}
(r_{m_1}, r_{m_2}),r_{m_3}) \cdots) , r_{m_{j-1}}, r_{m_j}).
\] 
In particular, 
\begin{equation}\label{w1iszero}
w_1=r_1-\frac 12 r_0^2 = 0
\end{equation}
 and 
\[\aligned
w_2  &=\frac{i}{2\pi}\left(r_2 - r_0r_1 + \frac 13 r_0^3\right) \\
& = -\frac{1}{4\pi} \int_0^{2\pi}\int_0^{2\pi -s} \PB{V(\phi_s(x,p))}{V(\phi_{s+u}(x,p))}\  
duds \\
& =-\frac{1}{4\pi} \int_0^{2\pi} \int_0^u \PB{V(\phi_s(x,p))}{V(\phi_u(x,p))}\ dsdu.
\endaligned\]
In what follows we will let, for any function $F$ on phase space,
\begin{equation}\label{2ndAve}
F^{\Delta}(x,p) :=-\frac{1}{4\pi} \int_0^{2\pi} \int_0^u \PB{F(\phi_s(x,p))}{F(\phi_u(x,p))}\ 
dsdu,
\end{equation}
so that 
\begin{equation}\label{whatw2is}
w_2 = V^{\Delta}.
\end{equation}

We are now in a position to extend Theorem \ref{BandInv1Thm} and 
obtain higher order terms in the asymptotic expansion of 
\begin{equation}\label{thetrace}
  \sum_j f(E_j) \sum_{k=1}^{m_j} \varphi(\mu_{k,j}) = \tr\left[ f(S_\chi)\varphi(W) \right].
\end{equation} 

\begin{theorem}\label{theoremX}
 Let $f\in C_{0}^{\infty}(\bbR)$ and  $\varphi(s) = s^{l+1}$ in 
(\ref{thetrace}), where $l\geq 0$ is an integer.  
Then (\ref{thetrace}) has a semiclassical expansion in powers of $\h$, and 
the second coefficient is equal to 
\[
(l+1)\int f(\frac{|x|^2+|p|^2}{2})(V^{\av})^l V^{\Delta}\ dxdp +Q^f_l ,
\]
where $Q^f_l $ is equal to the integral over $\bbR^{2n}$ of the following expression:
\begin{multline}\label{secondinvariant}
B_{2}(f(H_{0}), (V^{\av})^{l+1}) + f(H_{0})(V^{\av})^{[l]} + \\ 
+(V^{\av})^{l+1}\left[ f'(H_{0})(V+V^{\av}) -\calR(f)(H_{0})\right].
\end{multline}
Here $B_{2}$ is the $\h^2$ operator in the Moyal product, we have let
\[
w^{[l]} = \sum_{j=0}^{l-1}w^{j}B_{2}(w, w^{l-j}),
\]
and $\calR$ is the operator
\[
\calR(f)(s) = \frac{n}{8}f''(s) + \frac{s}{12}f'''(s).
\]
\end{theorem}
The proof is a calculation that we have sketched in the appendix.

\subsection{The case of odd potentials}

Recall that $V$ is odd iff $V^{\av}$ is identically zero.  By equations (\ref{w1iszero},
\ref{whatw2is}), in this case the operator $W$ is of order $-2$ and its principal symbol
is $V^{\Delta}$.  
\begin{theorem}\label{OddBandInvariant}  For odd potentials $V$, the integrals
\[
\int f(\frac{|x|^2+|p|^2}{2})\,\varphi(V^{\Delta})\ dxdp
\]
where $f\in C_{0}^{\infty}(\bbR)$ and $\varphi\in C^{\infty}(\bbR)$
are spectral invariants.
\end{theorem} 
\begin{proof}
Analogously as in Theorem \ref{BandInv1Thm}, the previous quantity
is the coefficient of the leading order term in the asymptotic expansion of
  \[ 
  \sum_j f(E_j) \sum_{k=1}^{m_j} \varphi(\h^{-2}\,\mu_{k,j}) 
  = \tr\left[ f(S_\chi)\varphi(\h^{-2}W) \right] . 
  \]
In the present case $\h^{-2}W$ is a pseudodifferential operator of order zero
and principal symbol $V^{\Delta}$.
\end{proof}
\subsection{Perturbations by semiclassical potentials} $ $

More generally, we can consider the harmonic oscillator perturbed by a semiclassical potential, 
that is, consider
\begin{equation}\label{semPotPer}
S = S_0+ \hbar^2 V(x, \hbar),
\end{equation}
where 
\[
V(x, \hbar) \sim V_0(x) + \hbar V_1(x) + \hbar^2 V_2(x) + \cdots
\]
as symbols in ${\mathcal O}(\langle x,\xi\rangle^{m})$.
The first band variant was calculated in \S 3.1, with $a(x, p)=V_0(x)$. For the higher invariants,  
a similar calculation as in \S 3.3 shows that for $k \ge 1$, 
\[
\dot r_k = -\{ H_0, r_k\} - i V_k - i\sum_{m+n+l+1=k} B_l(V_m, r_n)  \quad\text{and}\quad 
r_k|_{t=0} = 0.
\]
Again these equations can be solved via Duhamel's principle. In particular, we get 
\[
r_1(t, x, p) = \frac 12 r_0^2 + i \int_1^t V_1(\phi_s(x, p)) ds,
\]
which leads to $w_1=V_1^{\av}$.  More generally, 
\[
r_k(t, x, p) = i \int_0^t V_k(\phi_s(x,p)) ds + \mbox{terms depending only on\ }
V_0, \cdots, V_{k-1}
\]
and 
\[
w_k =V_k^{\av} +  \mbox{terms depending only on\ }V_0, \cdots, V_{k-1}.
\]
So if we take $\varphi(x)=x^{l+1}$, then the  $\hbar^{2l+k}$ term in 
  \[ \sum_j f(E_j) \sum_{k=1}^{m_j} \varphi(\mu_{k,j}(\h)) = \tr\left[ f(S_\chi)\varphi(W), 
\right]\]
is 
\[f(H_0) (V_0^{\av})^l V_k^{\av}  + \mbox{terms depending only on\ }V_0, \cdots, V_{k-1}.\]
We conclude:
\begin{theorem}\label{1stinvsemipotentials}
From the $k^{th}$ term of the expansion one can spectrally determine the quantities
\begin{equation} \label{kth}
\int f(\frac{|x|^2+|p|^2}2) (V_0^{\av})^l V_k^{\av} dx dp + Q^f_l(V_0, \cdots, V_{k-1}),
\end{equation}
where $Q^f_l(V_0, \cdots, V_{k-1})$ depends only on $V_0, \cdots, V_{k-1}$ (and their 
derivatives) and on $f$, $l$. 
\end{theorem}

\section{First inverse spectral results}

\subsection{Spectral rigidity} $ $

Obviously for any $O \in SO(n)$, the rotated potential 
\[V^O(x):=V(Ox)\] 
is isospectral with the potential $V(x)$. From this observation one can easily construct 
trivial families of isospectral potentials. 
\begin{definition}
We say a potential $V(x)$  
is \emph{spectrally rigid} if for any smooth family of isospectral potentials $V^t(x)$,  
with $V^0(x)=V(x)$, 
there is a smooth family of orthogonal matrices $O_t \in SO(n)$,  such that 
$V^t(x) = V(O_t x )$. 
\end{definition}

\begin{remark}
Consider the Taylor expansion of $V$ at the origin:
\[
V(x)  \sim V(0) + \sum \frac{\partial V}{\partial x_i}(0) x_i + \sum \frac{\partial^2 V}{\partial 
x_i\partial x_j}(0) x_i x_j +  \mbox{higher order terms}.
\]
It is obvious that if we rotate $V(x)$ to $V^O(x)$, the following remain unchanged:
\begin{itemize}
\item the constant term $V(0)$, 
\item the length $|\nabla V(0)|^2 = \sum (\frac{\partial V}{\partial x_i}(0))^2$,
\item the eigenvalues of the quadratic form 
$\sum \frac{\partial^2 V}{\partial x_i\partial x_j}(0)x_{i}x_{j}$.
\end{itemize}
We will see that these data are spectrally determined in many cases.
Moreover, it is easy to see that
 the effect that the constant term $V(0)$ and the linear term 
 $\sum \frac{\partial V}{\partial x_i}(0) x_i$ have on the spectrum are simply translations, 
 and therefore one may drop these terms and only consider potentials of the form 
\[V(x) = \sum a_i x_i^2 + \mbox{higher order terms}.\] 
(For a semiclassical potential, one may assume the leading term $V_0(x)$ if of this form. )
\end{remark}
\begin{remark}
For a semiclassical potential one has another method to construct trivial isospectral families. 
Namely, suppose 
\[
V(x, \hbar) = V_0(x)+ \hbar V_1(x)+\hbar^2 V_2(x) + \cdots 
\]
is a semiclassical potential. Then a change of variables 
\[ 
x_i \to x_i + b_i \hbar^2 
\]
converts the operator 
\[
\frac 12( -\hbar^2 \Delta + |x|^2) + \hbar^2 V(x, \hbar)
\] 
into 
\[
\frac 12( -\hbar^2 \Delta + |x|^2) + \hbar^2 \left(V(x+b\hbar^2, \hbar) + 
\sum x_i b_i + \hbar^2 \frac{b^2}2\right).
\] 
Therefore the semiclassical potential 
\[
\widetilde V(x, \hbar) = V(x+b\hbar^2, \hbar) + \sum x_i b_i + \hbar^2 \frac{b^2}2
\]
is isospectral with $V(x, \hbar)$. Note that if we take $b_i = -\frac{\partial V_0}{\partial x_i}
(0)$, then the leading term of $\widetilde V(x, \hbar)$ is 
\[\widetilde V_0 (x) = V_0(x) - \sum \frac{\partial V_0}{\partial x_i}(0) x_i,\]
which has no linear term. Iteratively using this method, i.e. making changes of variables 
\[x_i \to x_i + b_i x^m,\]
$m \ge 2$, one can convert $V(x, \hbar)$ to an isospectral semiclassical potential
\[\widehat V(x, \hbar) = \widehat V_0(x) + \hbar \widehat V_1(x) + \hbar^2 \widehat V_2(x) + 
\cdots \]
where each $\widehat V_i(x)$ has no linear term.  
\end{remark}

Now let $V(x)$ be an even potential function. We say that $V$ is \emph{formally spectrally 
rigid} if for any smooth family of isospectral even potentials,  $V^t(x)$, with $V^0(x)=V(x)$, 
we have 
\begin{equation}\label{FSR}
\left.\frac {d^k}{dt^k}\right|_{t=0} V^t= 0.
\end{equation}
for all $k \ge 1$. 
\begin{theorem}
Let $V(x)$ be a potential which is even in each variable. 
Then $V$ is formally spectrally rigid if it satisfies: 
the only even function $W$ such that 
\begin{equation}\label{FSR2} 
\int f(|z|^2) \varphi(V^{\av}) W^{\av} dzd\bar z = 0
\end{equation}
holds for all $f \in C_0^\infty(\mathbb R)$ and $\varphi \in C^\infty(\mathbb R)$  is the zero 
function, $W = 0$. 
\end{theorem}
\begin{proof}
Let $V^t(x)$ be a smooth family of isospectral even potentials, such that $V^0(x)=V(x)$. Then
\[
\aligned
0 = & \left.\frac d{dt}\right|_{t=0} \int f(|z|^2) \varphi\left((V^t)^{\av}\right) dzd\bar z \\
& = \int f(|z|^2) \varphi'\left(V^{\av}\right) \left.\frac d{dt}\right|_{t=0} (V^t)^{\av} dz d\bar z,
\endaligned
\] 
which implies $\frac d{dt}|_{t=0}V^t = 0$. Similarly 
\[
\aligned
0 = & \left.\frac {d^2}{dt^2}\right|_{t=0} \int f(|z|^2) \varphi\left((V^t)^{\av}\right) dzd\bar z \\
& = \int f(|z|^2) \left[\varphi''\left(V^{\av}\right)\left( \left.\frac d{dt}\right|_{t=0} (V^t)^{\av} 
\right)^2
 + \varphi'(V^{\av}) \left.\frac{d^2}{dt^2} \right|_{t=0} (V^t)^{\av} 
 \right] dzd\bar z,  
\endaligned
\] 
which implies (\ref{FSR}) for $k=2$. 
Continuing in the same fashion, one can prove (\ref{FSR}) for all $k$.
\end{proof} 

Similarly one can define the \emph{formally spectral rigidity} for semiclassical potentials
\[V(x, \hbar) = V_0(x) + \hbar V_1(x) + \hbar^2 V_2(x) + \cdots,\]
where each $V_k(x)$ is even, by requiring that for any smooth familhy of  isospectral 
semiclassical potentials $V^t(x, \hbar)$, one has
\[
\left.\frac {d^k}{dt^k}\right|_{t=0} V_j^t= 0
\]
for all $k \ge 1$ and all $j \ge 0$. And one can show by using Theorem 3.4 that $V(x, \hbar)$ is 
formally spectral rigid if 
\begin{enumerate}
\item $V_0$ is formally spectrally rigid, e.g. it satisfies the hypotheses of Theorem 4.4 above. 
\item For each $j \ge 1$, $V_j$ satisfies the following condition: the only even function $W$ 
such that 
\[
\int f(|z|^2) (V^{\av})^l W^{\av} dzd\bar z = 0
\]
holds for all $f \in C_0^\infty(\mathbb R)$ and $l \in \mathbb N$  is the zero function. 
\end{enumerate}

\subsection{Recovering one dimensional even potentials} $ $

According to (\ref{AveMon}),  if $V$ is an analytic function of one variable, 
then after averaging the only terms that survive are the even terms in the Taylor expansion of 
$V$.  It follows that if $n=1$ the semiclassical spectrum of $S$ determines the even part of $V$ 
for analytic functions $V$.  The following theorem shows that the analyticity assumption can
be dropped:

\begin{theorem}\label{OneDimEven}
For $n=1$, the first band invariant determines the even part of $V$, i.e. $V(x)+V(-x)$, for any 
smooth perturbation $V$. 
\end{theorem}
\begin{proof}
It is sufficient to prove the statement for $V$ even. 
Then from the first band invariant gives us, for any $r>0$, the integral
\[\int_0^{\pi/2} V(r\cos \theta)\ d\theta.\]
Making the change of variables $s=r\cos \theta$ and $u = s^2$, and denoting 
$V_1(x) = V(\sqrt x)/\sqrt{x}$, the above integral becomes 
\[
\int_0^r V(s)\frac 1{\sqrt{r^2-s^2}}\ ds = \int_0^{r^2} V(\sqrt u) \frac 1{\sqrt{r^2-u}}\ 
\frac{du}{\sqrt{u}} =  \int_0^{r^2} V_1(u) (r^2-u)^{-1/2}\ du.
\]
The latter one equals 
\[\Gamma(\frac 12) (J^{\frac 12}V_1)(r^2),\]
where $J^{\frac 12}$ is the fractional derivative of order $\frac 12$. So if we apply $J^{\frac 
12}$ again and integrate, we can recover $V_1$, and thus $V$ itself.  
\end{proof}

\begin{remark}
In \cite{GUW}, we showed that by using semiclassical invariants modulo $O(\hbar^4)$, one 
can spectrally determine not only $V(x)+V(-x)$, but also $V^2(x)+V^2(-x)$. It follows that can 
determine $V(x)$ itself under suitable symmetry conditions. It turns out the same result holds in 
higher dimensions, i.e. one can spectrally determine the integral $\int V(x) d\sigma_x$ as well 
as the integral $\int V^2(x)d\sigma_x$ over any sphere $|x|=r$.  In particular, one can 
distinguish radially symmetric potentials from other potentials. For more details, c.f. \cite
{GUW}.
\end{remark}

\subsection{Recovering one dimensional odd analytic potentials} 

Now assume $V$ is a one dimensional odd analytic potential,
\[
V(x) = a_1 x + a_3 x^3 + \cdots.
\]
Then
\[
V^{\Delta} (z) =- \frac 1{4\pi} \int_0^{2\pi} \int_0^u \left\{\sum_{k \mbox{\scriptsize\ odd}} 
a_k (\frac{e^{is}z+e^{-is}\bar z}2)^k\,,\,\sum_{l \mbox{\scriptsize\ odd}} a_l (\frac{e^{iu}z
+e^{-iu}\bar z}2)^l\right\}\; dsdu.
\]
From the first band invariant for odd potentials (see Theorem \ref{OddBandInvariant}),
\[
\int e^{-\mu |z|^2} V^{\Delta} dz d\bar z,
\]
we get, by replacing $\mu$ by $\frac{\mu}{\lambda^2}$ and making change of variables $z  \to 
\lambda z$, the spectrally determined expression
\begin{equation}\label{xkxl}
\sum_{k, l \mbox{\scriptsize\ odd}} a_k a_l \lambda^{k+l+2} \int_{\bbC} e^{-\mu |z|
^2}\int_0^{2\pi} \int_0^u \PB{( {e^{is}z+e^{-is}\bar z} )^k}{( {e^{iu}z+e^{-iu}\bar z} )^l}  
dsdu\ dzd\bar z.
\end{equation}
Note that in complex coordinates,
\[
\PB{f}{g} = \Im \frac{\partial f}{\partial z}\frac{\partial g}{\partial \bar z},
\]
so we have for $k=1$ and $l=2l_1+1$ odd,
\[\aligned
& \int_{\bbC} e^{-\mu |z|^2}\int_0^{2\pi} \int_0^u \PB{ {e^{is}z+e^{-is}\bar z} }{( {e^{iu}z
+e^{-iu}\bar z} )^l}  dsdu\ dzd\bar z \\
= &  \int_{\bbC} e^{-\mu |z|^2}\int_0^{2\pi} \int_0^u \Im l e^{i(s-u)} ( {e^{iu}z+e^{-iu}\bar 
z} )^{l-1}dsdu\ dzd\bar z \\
= & l {l-1 \choose (l-1)/2}    \int_{\bbC} e^{-\mu |z|^2} |z|^{l-1}dzd\bar z \int_0^{2\pi} 
\int_0^u \sin (s-u) dsdu \\
= &- 2\pi^2  l {l-1 \choose l_1}  \frac{l_1!}{\mu^{l_1+1}}\\
= &- \frac{ 2\pi^2 l!}{l_1! \mu^{l_1+1}}.
\endaligned\]
Similarly for $k=2k_1+1$ and $l=1$ we have
\[
\int_{\bbC} e^{-\mu |z|^2}\int_0^{2\pi} \int_0^u \PB{ {(e^{is}z+e^{-is}\bar z)^k} }{ {e^{iu}z
+e^{-iu}\bar z} }  dsdu\ dzd\bar z  = - \frac{2\pi^2 k!}{k_1! \mu^{k_1+1}}.
\]
In particular, we get from the lowest order term of $\lambda$ in (\ref{xkxl}) (i.e.\  the 
coefficient of $\lambda^4$) the number $\frac{2\pi^2}{\mu} a_1^2$ from which one can 
recover $a_1$ up to a sign.  At this point we must distinguish two cases.

Case 1: $a_1 \ne 0$. Suppose we have recovered the first $m$ coefficients $a_1, \cdots, a_
{2m-1} $,  up to an overall sign (the same for all). Then by looking at the coefficient of $
\lambda^{2m+4}$ in the expression (\ref{xkxl}) one can recover the number $a_1 a_{2m+1} + 
a_{2m+1} a_1$, since all other products of pairs are known precisely. 
It follows that $a_{2m+1}$ is determined up to the same sign ambiguity as 
$a_1$. 

Case 2:  $a_1 =0$. More generally, suppose we have already found that $a_1=\cdots=a_{k-2}=0$, 
where $k \ge 3$ is an odd number. Let $l \ge k$ be odd too. Then  a calculation as above shows 
that 
\[\aligned
& \int_{\bbC} e^{-\mu |z|^2}\int_0^{2\pi} \int_0^u \PB{ (e^{is}z+e^{-is}\bar z)^k }{( {e^{iu}
z+e^{-iu}\bar z} )^l}  dsdu\ dzd\bar z \\
= & C \sum_{m=0}^{k-1} \frac 1{k-2m} {k-1 \choose k-1-m} {l-1 \choose \frac{l-k+2m}2},
\endaligned\]
where $C$ is some non-vanishing constant that depends on $\mu$.  Note that 
\[
\frac 1k  {k-1 \choose k-1} {l-1 \choose \frac{l-k}2} > 0
\]
and for $0<m<\frac k2$,
\[
 \frac 1{k-2m} {k-1 \choose k-1-m} {l-1 \choose \frac{l-k+2m}2}+\frac 1{2m-k} {k-1 \choose 
m-1} {l-1 \choose \frac{l+k-2m}2}>0, 
\]
so the expression 
\[
 \int_{\bbC} e^{-\mu |z|^2}\int_0^{2\pi} \int_0^u \PB{ (e^{is}z+e^{-is}\bar z)^k }{( {e^{iu}z
+e^{-iu}\bar z} )^l}  dsdu\ dzd\bar z 
\]
is non-vanishing as well. It follows that for $l=k$, one can determine $a_k^2$, and for $l>k$, 
one can determine $a_k a_l$, from which one can determine $a_l$ up to the same sign 
ambiguity as $a_k$.

In conclusion, we have proved:
\begin{theorem}
For any odd analytic potential $V$, in one dimension,
the first band invariant for odd potentials determines $V(x)$ up to a sign.
\end{theorem}

\subsection{Recovering constant and quadratic terms} $ $

In all dimensions  $V(0)$ is spectrally determined, since 
\[V(0) = V^{\av}(0) = \lim_{r \to 0} \int_{|x|^2+|p|^2=r^2} V^{\av} d\lambda.\]

About recovering the quadratic terms, we have already seen that we can only hope to recover 
the eigenvalues of the Hessian at the origin.  Since the average of any odd function vanishes, 
we can assume that $V(x)$ is even itself. 
We will assume moreover that $V(x)$ is analytic, so that up to a rotation, 
\[
V(x) = V(0) + \sum a_i x_i^2 + \widetilde V(x),
\] 
where $\widetilde V(x)$ contains only terms that are homogeneous in
 $x$ of degree greater than four. It follows that
\[
V^{\av}(x) = V(0) + \sum a_i |z_i|^2 + \widetilde V^{\av}(z),
\]
where $\widetilde V^{\av}(z)$ contains only terms that are homogeneous in $z$ of degree 4
and higher.  So from the spectral invariant 
  \[ \int e^{-{\mu} |z|^2} V_{\av}^{k+1}\, dzd\bar z \]
we get, by replacing $\mu$ by $\frac{\mu}{\lambda^2}$ and making change of variables $z_i 
\to \lambda z_i$, the expressions
  \[ 
     \lambda^{2k+4} \int e^{-{\mu} |z|^2 }\left( \sum a_i |z_i|^2 \right)^{k+1}\, dzd\bar z   
   +  \mbox{terms involves at least\ } \lambda^{2k+6} .
   \]
It follows that the integral 
\[ \int e^{-{\mu} |z|^2 }\left( \sum a_i |z_i|^2 \right)^{k+1}\, dzd\bar z  \] 
is spectrally determined for all $k$. In particular, one can determine 
\[
\int e^{-\mu |z|^2} e^{\sum a_i |z_i|^2} dzd\bar z = \pi^n \frac{1}{\mu-a_1} \cdots \frac{1}
{\mu-a_n},
\]
and thus determine the polynomial 
\[
(\mu-a_1) \cdots (\mu-a_n).
\]
This proves:
\begin{theorem}
For analytic potentials the eigenvalues of the Hessian at the origin are
spectrally determined.
\end{theorem}

\subsection{Recovering the linear term for certain potentials}

First assume $V$ is an analytic odd potential, so that
\[
V(x) =V_1(x) + V_3(x) + \cdots,
\]
where $V_k$ is homogeneous of degree $k$, and in particular, $V_1(x) = \nabla V(0) \cdot x$.  Then according to Theorem 3.4, the integrals
\[
\int e^{-\mu |z|^2} V^{\Delta} dzd\bar z
\]
are spectrally determined. Again we replace $\mu$ by $\frac{\mu}{\lambda^2}$ and make change of variables $z \to \lambda z$, we get the expression
\[
\lambda^2 \int e^{-\mu |z|^2} V_1^{\Delta} dzd\bar z + \mbox{terms involves at least \ } \lambda^4
\]
So one can spectrally determine the integral
\[ \int e^{-\mu |z|^2} V_1^{\Delta} dzd\bar z \]
from which one can read off $\| \nabla V(0) \|^2$.

One can improve this result slightly by considering analytic potentials of the form
\begin{equation}\label{noV2V4}
V(x)=V_0(x)+V_1(x)+V_3(x)+V_5(x)+V_6(x)+ \cdots,
\end{equation}
i.e. analytic potentials with vanishing quadratic and quartic terms.  Note that this property,
\[
V_2=0 \quad \mbox{and} \quad V_4=0,
\]
is a spectrally property: We have already seen in section 4.4 that one can determine whether $V_2=0$ using spectral data. In this class of potentials, one
can spectrally determine the integrals
\[ \int _{S_r} (V_4^{\av}) d\lambda.\]
If all these integrals vanish, then we get
\[\int_{S_r} e^{itV_4^{\av}} d\lambda = \mathrm{Vol}(S_r),\]
from which one can deduce that $V_4^{\av}=0$, and thus $V_4=0$.

Now suppose $V(x)$ is a potential of the form (\ref{noV2V4}). Without loss of generality, one can assume that the constant term is zero. 
Then according to theorem 3.3, with $l=0$ and $f(t)=e^{-\mu t}$, one can spectrally determine
\[
\int \left[ e^{-\mu |z|^2} V^{\Delta} + B_2(e^{-\mu |z|^2}, V^{\av}) \right] dzd\bar z.
\]
Since $V_2^{\av}=V_4^{\av}=0$, after replacing $\mu$ by $\frac{\mu}{\lambda^2}$ and making change of variables $z \to \lambda z$, the lowest order term is again
\[ 
\lambda^2 \int e^{-\mu |z|^2} V_1^{\Delta} dzd\bar z,
 \]
and from this one gets $\| \nabla V(0) \|^2$.
In conclusion, we have
\begin{theorem}
Let $V(x)$ be an analytic potential. Then one can spectrally determine whether $V(x)$ is of the form (\ref{noV2V4}), and if $V$ is in this class, one can spectrally determine its linear term $V_1(x)$ up to rotation.
\end{theorem}

\section{The first band invariant in dimension 2}

 The averaging procedure on the space of functions
 \[
C^\infty_{\mathrm{even}}(\bbR^2) := \{  V(x_1, x_2) = 
\tilde{V}(x_1^2, x_2^2), \  \tilde{V}\in C^\infty(\bbR^2)\}
\]
 of smooth functions on $\bbR^2$ which are even in each variable has particularly nice
 properties that we investigate in this section.
 
Let $\pi$ be the cotangent fibration $\pi:\mathbb C^2 = T^* \mathbb 
R^2 \to \mathbb R^2$, and let $V^{\av}$ be the average of $\pi^* V$ with respect to the circle 
action 
  \begin{equation}\label{4.1}
  e^{i\theta} z  = (e^{i\theta}z_1, e^{i\theta} z_2).
  \end{equation} 
As we have seen, for $\varphi \in C^\infty(\mathbb R)$, the integrals 
  \begin{equation}\label{4.2}
  \varphi_\mu(V) = \int e^{-\frac {\mu}2 (|z_1|^2+|z_2|^2)} \varphi(V^{\av})\, dzd\bar z
  \end{equation}
are spectral invariants of the operator (\ref{GenPer}).

To analyze these invariants we'll begin by decomposing $V^{\av}$ into its Fourier coefficients 
with respect to the circle action 
  \begin{equation}\label{4.3}
  e^{i\theta}z = (e^{i\theta}z_1, e^{-i\theta}z_2),
  \end{equation}
and this we'll do by first examining the one-dimensioanl analogues of these Fourier coefficients. 
More explicitly, let 
  \begin{equation}\label{4.4}
  B_r: C^\infty_{\mathrm{even}}(\mathbb R) \to C^\infty_{\mathrm{even}}(\mathbb R)
  \end{equation}
be the operator 
\[
f(x) = f\left(\frac{z + \bar z}2\right)  =g(|z|^2 e^{2i\theta}) \mapsto g_r,
\]
where $g_r$ is the $2r$-th Fourier coefficent of $g(|z|^2 e^{2i\theta})$ with respect to $\theta$. 
By definition this operator commutes with the homothety $x \to \lambda x$, $\lambda \in 
\mathbb R_{+}$, and hence maps $x^{2k}$ into a multiple, $\gamma_{k,r} x^{2k}$, of itself. 
To compute $\gamma_{k,r}$ we note that 
\[
\left( \frac{z+\bar z}2 \right)^{2k}  = \frac 1{4^k} \sum_{r=-k}^{r=k} {2k \choose r+k} |z|^
{2k} e^{2i r\theta}.
\]
Hence 
  \begin{equation}\label{4.5}
  B_r x^{2k} = \frac 1{4^k} {2k \choose r+k} x^{2k}
  \end{equation}
for $|r| \le k$ and zero for $|r|>k$.  Note that $B_{0}$ is invertible.

To describe the asymptotic dependence of $\gamma_{k,r}$ 
on $k$ for $k >> 0$ we apply Stirling's formula 
  \[n! \sim \sqrt{2\pi n}\left(\frac ne\right)^n \left( 1 + a_1 n^{-1} + a_2 n^{-2} + \cdots \right) \]
to the quotient 
  \[
  \aligned 
  \frac 1{4^k} {2k \choose k+r} & = 4^{-k} \frac{(2k)!}{(k+r)!(k-r)!} \\
  & = \frac{4^{-k}\sqrt{4\pi k}(2k/e)^{2k}}{\sqrt{4\pi^2 (k^2-r^2)}(\frac{k+r}e)^{k+4} (\frac
{k-r}e)^{k-r}} \left( 1+ b_{1,r}k^{-1} + \cdots \right)\\
  &=  \frac 1{\sqrt{\pi k}} \left(  1+ c_{1,r}k^{-1} + \cdots \right),
  \endaligned
  \]
giving the asymptotic expansion
  \begin{equation}\label{4.6}
  \gamma_{k,r} \sim \frac 1{\sqrt{\pi k}} \sum_{l=0}^\infty c_{l,r} k^{-l}, \quad c_{0,r} = 1,
  \end{equation}
or (more intrinsically) the asymptotic expansion 
  \begin{equation}\label{4.7}
  B_r  = \sum_{r=0}^{\infty} c_{k,r} Q^{-r-\frac 12}
  \end{equation}
where $Q = \frac x2 \frac d{dx}$. Hence in particular $B_r$ is a pseudodifferential operator on 
$C^\infty_{\mathrm{even}}(\mathbb R)$ of order $-\frac 12$. 

\medskip
Coming back to the problem in two dimensions that prompted these computations, 
for each $V\in C^\infty_{\mathrm{even}}(\mathbb R^2)$ let
$R_nV$ the $n$th Fourier coefficient of $V^{\av}$ with respect to the circle action 
 (\ref{4.3}). We will prove:
\begin{theorem}
The operator
  \begin{equation}\label{4.7a}
 R_n: C^\infty_{\mathrm{even}}(\mathbb R^2) \to C^\infty_{\mathrm{even}} (\mathbb R^2)
 \end{equation}
is zero for $n \ne 0 \mod 4$, and modulo the identification $C^\infty_{\mathrm{even}} = 
C^\infty_{\mathrm{even}}(\mathbb R) \hat \otimes C^\infty_{\mathrm{even}}(\mathbb R)$,
  \begin{equation}\label{4.8}
  R_n  = B_r \hat \otimes B_r, \quad n=4r
  \end{equation}
(i.e. for functions of the form $f(x_1, x_2) = g(x_1) h(x_2), R_n f = B_r g(x_1) B_r h(x_2)$).  
In particular $R_{0}$ is invertible.
\end{theorem}
\begin{proof}
It suffices to check this for $g=x_1^{2k}$ and $h = x_2^{2l}$ in which case 
 \[\aligned
  V^{\av} & = \frac 1{2\pi} 4^{-(k+l)} \int (z_1 e^{i\theta} + \bar z_1 e^{-i\theta})^{2k}(z_2 
e^{i\theta} + \bar z_2 e^{-i\theta})^{2l}\, d\theta \\
  & = \frac 1{2\pi}4^{-k-l} \sum_{s,t} {2k \choose s}{2l \choose t} z_1^s \bar z_1^{2k-s} 
z_2^t \bar z_2^{2l-t}  \int e^{2i(s+t-k-l)\theta} d\theta \\
  & = 4^{-k-l} \sum_{s+t = k+l} {2k \choose s}{2l \choose t} z_1^s \bar z_1^{2k-s} z_2^t \bar 
z_2^{2l-t} \\
  & = \sum_{|r| \le \min(k,l)} V_{k,l,r},
  \endaligned
  \]
where, in polar coordinates $z_j = |z_j|e^{i\theta_j}$,
  \begin{equation}\label{4.9}
  V_{k,l,r} = 4^{-(k+l)} {2k \choose k+r}{2l \choose l-r} |z_1|^{2k}  |z_2|^{2r} e^{2ir(\theta_1 
-\theta_2)}.
  \end{equation}
Hence setting $x_j = |z_j|$ we get for the $n$-th Fourier coefficient, $n = 4r$, of $V^{\av}$ 
with respect to the circle action (\ref{4.3}), 
\[(\gamma_{k,r} x_1^k)  (\gamma_{l,r}x_2^l) = B_r x_1^k B_r x_2^l.\] 
\end{proof}
\begin{remark}
The operator $A_r = R_n, n=4r$, can also be described as a Radon transform. Namely consider 
the double fibration 
\[
\xymatrix{
\; & \bbC\times\bbC  \ar[ld]_\pi \ar[rd]^\rho \\  
\bbR \times \bbR & & \bbR \times \bbR  }
\]
where $\pi(z_1, z_2) = (\mathrm{Re}z_1, \mathrm{Re} z_2)$ and $\rho(z_1, z_2) = (|z_1|, |
z_2|)$. Then $A_r: C^\infty_{\mathrm{even}} (\mathbb R^2) \to C^\infty_{\mathrm{even}}
(\mathbb R^2)$ is the transform 
  \begin{equation}\label{4.10}
  A_r V =  \rho_*\left( \int \tau^*_{\theta_1, \theta_2} \pi^* V e^{-4r\theta_2}\, d\theta_1 d
\theta_2\right),
  \end{equation}
where $\tau_{\theta_1, \theta_2}(z_1, z_2) = e^{i\theta_1} (e^{i\theta_2}z_1, e^{-i\theta_2}
z_2).$
\end{remark}

If we take $\varphi(t) = t^m$ in (\ref{4.1}) we get the spectral invariants
  \begin{equation}\label{4.8a}
  \int e^{-\frac{\mu}2 (|z_1|^2 +|z_2|^2)} (V^{\av})^m \, dzd\bar z
  \end{equation}
and if we expand $V^{\av}$ in its Fourier series with respect to the circle action (\ref{4.3}), i.e. 
express $V^{\av}$ as the sum 
  \[ \sum A_r V(x_1, x_2) e^{4ri\theta}, \]
then (\ref{4.8a}) becomes 
  \begin{equation}\label{4.9a}
  \int e^{-\frac{\mu}2 (x_1^2+ x_2^2)} \left( \sum_{r_1+ \cdots + r_m=0} A_{r_1}V \cdots A_
{r_m}V \right)\, dx_1 dx_2.
  \end{equation}
In particular for $m=1$ we obtain 
  \begin{equation}\label{4.10a}
  \int e^{-\frac{\mu}2 (x_1^2+ x_2^2)} A_0 V \, dx_1 dx_2,
  \end{equation}
and for $m=2$
  \begin{equation}\label{4.11}
  \int e^{-\frac{\mu}2 (x_1^2 + x_2^2)} \sum_{r} |A_r V|^2\, dx_1dx_2.
  \end{equation}

\section{Inverse spectral results in two dimensions} 

Throughout this section we work in two dimensions, and consider the inverse spectral problem 
for a perturbation of the harmonic oscillator by a potential (at times semicassical).

\subsection{Results for smooth perturbations}

\subsubsection{Potentials of the form $f_{1}(x_1^2)+f_{2}(x_2^2)$}

In this subsection we consider potentials of the form
\begin{equation}\label{formofthepotentials}
V(x_{1},x_{2})=f_{1}(x_1^2)+f_{2}(x_2^2).
\end{equation}
It is easy to see that
\begin{equation}\label{}
V^{\av} (x, p) =\varphi_{1}(x_{1}^{2}+p_{1}^{2 }) + 
\varphi_{2}(x_{2}^{2}+p_{2}^{2 }), 
\end{equation}
where
\[
\varphi_{j}(r) = \frac{\Gamma(1/2)}{\pi}\; J^{1/2}\left(\frac{f_{j}(s)}{\sqrt{s}}\right)(r)
\]
(see the proof of Theorem \ref{OneDimEven}).

\begin{theorem}
Generically, potentials of the form (\ref{formofthepotentials}) are spectrally determined,
up to the obvious symmetries: Exchanging of the roles of $x_1$ and $x_2$, and adding a
constant to $f_1$ and subtracting the same constant from $f_2$.
(The genericity condition is (\ref{themonotoneassumption})).
\end{theorem}
\begin{proof}
 
As we have seen (c.f.\ (\ref{Sz})), the integrals
  \begin{equation}\label{usedtobesixpointtwo}
  \int_{S_{r}} (V^{\av})^k(\theta)\, d\theta,
  \end{equation}
where $S_{r}\subset\bbR^{4}$ is the sphere of radius r,
are spectrally determined. 
We claim that the previous integral equals
\begin{equation}\label{whattheintegralequals}
2\pi^{2}r^{3}\int_{0}^{1}\left[ \varphi_{1}(r^{2}(1-u)) +\varphi_{2}(r^{2}u)\right]^{k}\,du.
\end{equation}
To see this, consider the map $\Phi:\bbR^{4}\to\bbR_{+}^{2}$ given by
\[
\Phi(x,p) = (x_{1}^{2}+p_{1}^{2}\,,\,x_{2}^{2}+p_{2}^{2}) =: (r^{2}u_{1}, r^{2} u_{2}).
\]
This maps the sphere $S_{r}$ onto the line segment $u_{1}+u_{2}=1$ (in the first
quadrant).  We will take $u=u_{2}$ as a coordinate in that segment. The fibers of $\Phi$ 
 are tori, which we can parametrize in complex 
coordinates by $(r\sqrt{u_{1}} e^{is}, r\sqrt{u_{2}} e^{it})$.  The variables $(s, t, u)$ 
parametrize $S_{r}$.  Using these variables to compute (\ref{usedtobesixpointtwo})
yields (\ref{whattheintegralequals}) (we'll omit the details of the calculation).

The integrals (\ref{whattheintegralequals}) determine the distribution function of
\[
\psi_{r}(u) = \varphi_{1}(r^{2}(1-u)) +\varphi_{2}(r^{2}u)
\]
on $[0,1]$.  We now make the genericity assumption that
\begin{equation}\label{themonotoneassumption}
\psi_{r}(u) \ \text{is monotone for each } r>0.
\end{equation}
and therefore its distribution function determines it up to the ambiguity $\psi_{r}(1-u)$, which 
amounts to switching the roles of $x_1$ and $x_2$.  Finally, it is not hard to see that knowing
the two-variable function $\psi$ determines $\varphi_1$ and $\varphi_2$ up to the
ambiguity of adding and subtracting a constant.
\end{proof}

\subsubsection{Semiclassical Potentials with quadratic $V_0$}

\begin{theorem}\label{SemiPots}
Semiclassical potentials of the form
\[
V(x, \hbar) = a x_1^2 + b x_2^2 + \hbar V_1 + \hbar^2 V_2 + \cdots,
\]
where $a\not =b$ and
$V_i\in C^\infty_{\mathrm{even}}(\mathbb R^2)$ for each $i$, are spectrally determined.
\end{theorem}
\begin{proof}
According to Theorem \ref{1stinvsemipotentials}, one can spectrally determine 
\[
\int e^{-\lambda(|z_1|^2+|z_2|^2) } \left( a|z_1|^2+b|z_2|^2\right)^l V_1^{\av} \ dzd\bar z
\]
plus a term not depending on $V_i$ with $i\geq 1$. We can rewrite this invariant as 
\[
\int e^{-\lambda (x_1^2+x_2^2)} \left(ax_1^2+bx_2^2\right)^l A_0V_1 dx_1 dx_2.
\]

From this we can determine
\[ 
\int e^{-\lambda (x_1^2+x_2^2)} e^{-\mu (ax_1^2+bx_2^2)} A_0 V_1 dx_1 dx_2 
\] 
for all  $\mu$ and $\lambda >0$ sufficiently large with respect to $\mu$
(so that the integral is absolutely convergent).   Note that the phase equals
\[
-(\lambda +a\mu)x_{1}^{2} - (\lambda + b\mu)x_{2}^{2},
\]
and since $a\not= b$ we can make a linear change of variables, 
$(\lambda, \mu)\to (\mu_{1}, \mu_{2})$ so that the invariant becomes
\[
\int e^{-\mu_{1}x_{1}^{2} - \mu_{2}x_{2}^{2}}\; A_{0}V_{1}\, dx_{1}\, dx_{2}.
\]
Integrating in polar coordinates
we get, up to a universal constant, the spectral invariant
\[ 
\int e^{-\mu_1 r_1 -\mu_2 r_2} A_0 V_1(\sqrt r_1, \sqrt r_2)/(\sqrt r_1 \sqrt r_2) dr_1 dr_2.
 \] 
This function of $\mu_{1}$ and $\mu_{2}$ is known spectrally in some infinite wedge
of $\bbR^{2}_{+}$ (determined by the condition that $\lambda >0$).  However, it
 is easy to see that this is an analytic function of $(\mu_{1}, \mu_{2})$, and therefore
 it is spectrally determined in the whole quadrant $\bbR^{2}_{+}$.
Using the inverse Laplace transform, one can determine $A_0 V_1 (x_1, x_2)$ pointwise, 
and thus determine $V_1$ itself. 

By an inductive argument, one can similarly spectrally determine each of the remaining
$V_k$. 
\end{proof}

\subsection{Inverse Spectral Results for Real Analytic Perturbations}

\subsubsection{Spectral rigidity of the quadratic potentials $ax_1^2 + bx_2^2$.}
In this section we will prove the following result, which will be used later:
\begin{proposition}
If $a \ne b$ the potential $V = ax_1^2 + bx_2^2$ is
``formally" spectrally rigid; i.e. if $V_t, -\varepsilon < t < \varepsilon$, is a smooth family of 
potentials each of which is an even function of $x_1$ and $x_2$ and has the same spectrum as 
$V$, then 
  \begin{equation}\label{7.1}
  \left.\frac{d^k}{dt^k}V_t\right|_{t=0} =0
  \end{equation}
for all $k$.
\end{proposition} 
\begin{proof}
We first note that if (\ref{7.1}) is non-zero for some $k$ then we can 
reparametrize $t$ so that (\ref{7.1}) is non-zero for $k=1$. Thus it suffices to prove that the 
function $W  = \frac{dV}{dt}|_{t=0}$ is zero. To see this we first note that by remark 5.2, $A_i V=0$ for $i \ne 0$ and $A_0V = V$. Thus by inserting $V_t$ into (\ref{4.9a}), 
differentiating with respect to $t$ and setting $t=0$, we get 
  \[ \int e^{-\frac{\mu}2 (x_1^2 + x_2^2)} V^{m-1}A_0 W\, dx_1dx_2 = 0 \]
for all $m$. Thus 
  \[ \int e^{-\frac{\mu}2 (x_1^2+x_2^2)} e^{-\frac{\nu}2 V}A_0 W \, dx_1dx_2 = 0 \]
and since 
  \[ -\frac{\mu}2 (x_1^2+x_2^2) - \frac{\nu}2 V = -\frac{\mu+\nu a}2 x_1^2 - \frac{\mu + \nu 
b}2 x_2^2 \]
and $a \ne b$, we can, by making a linear change of coordinates, $(\mu, \nu) \to (\mu_1, 
\mu_2)$, convert this equation into the form 
  \[ \int e^{-\frac{\mu_1}2x_1^2 -\frac{\mu_2}2 x_2^2} A_0 W(x_1^2, x_2^2)\, dx_1dx_2 = 
0\]
for all $\mu_1, \mu_2 $ in an infinite wedge of the 1st quadrant.  Arguing just as in the
end of the proof of Theorem \ref{SemiPots} we can apply the inverse Laplace transform
and conclude that $A_0W = 0$, and by the injectivity of $A_0$, that $W=0$. 
\end{proof}

\subsubsection{Spectral determinacy}

Since one can spectrally 
determine the integrals of $\sum |A_iV|^2$ over the circles $x_1^2+x_2^2 = t$, in particular 
one can determine whether this sum vanishes to infinite order at $t=x_1=x_2=0$, 
and hence whether 
$V$ itself vanishes to infinite order at $t=0$.   This can be strengthened, as follows:

\begin{theorem}
If $V$ is an analytic function on $\bbR^{2}$ such that its Taylor expansion at the origin 
is of the form
  \[
  V=V_2 + V_4 + V_6 + \cdots
  \]
with each $V_j$ homogeneous of degree $j$, and the quadratic term $V_2$ has two distinct 
eigenvalues $a \ne b$, then $V$ is spectrally determined (up to a rotation)\footnote
{Results of this nature for the Schr\"odinger operator, $-\h^2 \Delta_{\mathbb R^n}+V$, $V$ 
real analytic, can be found in \cite{GU}, \cite{GPU} and \cite{Hez}.  However, these 
results require strong ``non-rationality" assumptions on the coefficients of the leading term $
\sum a_ix_i^2$ of $V$. }. 
\end{theorem}
\begin{proof}
We can assume without loss of generality  that $V_2 = a x_1^2+b x_2^2$,
and let
  \[ 
  V^{\av} = a|z_1|^2 + b|z_2|^2 + V_4^{\av} + V_6^{\av} + \cdots 
  \]
be the Taylor series expansion of $V^{\av}$. From the spectral invariants 
  \[ \int e^{-\frac{\mu}2(|z_1|^2+|z_2|^2)} (V^{\av})^{k+1}(z_1, z_2)\, dzd\bar z \]
we get, by replacing $\mu$ by $\frac{\mu}{\lambda^2}$ and making change of variables $z_i 
\to \lambda z_i$, the expressions
  \[\aligned
   & \lambda^{2k+4} \int e^{-\frac{\mu}2(|z_1|^2+|z_2|^2)}\left( a|z_1|^2 + b|z_2|^2 \right)^{k
+1}\, dzd\bar z  \\
  & \quad +  (k+1)\lambda^{2k+6} \int e^{-\frac{\mu}2(|z_1|^2+|z_2|^2)}\left( a|z_1|^2 + b|z_2|
^2 \right)^k V_4^{\av}\, dzd\bar z
   + O(\lambda^{2k+8})
  \endaligned\]
and hence, from the proof of the spectral rigidity  of $ax_1^2 + bx_2^2$ and linearity 
of this expression in $V_{4}^{\av}$, this suffices to determine $V_4^{\av}$. 

Continuing: the next unknown term in the expansion above is
  \[ (k+1)\lambda^{2k+8} \int e^{-\frac{\mu}2 (|z_1|^2+|z_2|^2) } \left( a|z_1|^2 + b|z_2|^2 
\right)^k V^{\av}_6 \, dzd\bar z \]
and hence, by the same argument, this determines $V_6^{\av}$ (and 
by induction and repetition of this argument $V_8^{\av}, V_{10}^{\av}, \cdots$). 
Thus the Taylor series of $V^{\av}$ is spectrally determinable and hence by (\ref{AveMon}) 
so is the Taylor series of $V$. 
 \end{proof}
 
\subsubsection{The case of semiclassical potentials} 

Now consider perturbations of $S_0$ by the semiclassical potentials, 
\[ -\frac 12 (\hbar^2 \Delta + |x|^2) + \h^2(V_0 + \hbar V_1 + \hbar^2 V_2 + \cdots),\]
where 
\begin{equation}\label{vnot}
V_0 = a x_1^2 + bx_2^2 +\mbox{\ higher order terms\ }
\end{equation}
with $a \ne b$ as above, and all $V_i$'s are even functions. As we have seen from the previous 
section that the first spectral invariant is enough to determine $V_0$. To determine $V_i$'s with 
$i \ge 1$, we use the higher invariants.  We have seen in section 3.5 that the integral 
\[ 
\int e^{-\lambda (|z_1|^2 +|z_2|^2)} (V_0^{\av})^k V_1^{\av}\,dzd\zbar
\]
is spectrally determined. If we write
\[V_j = V^j_0 +V^j_2 +V^j_4 + \cdots\]
with each $V^j_i$ homogeneous of degree $i$, then the same argument as 
in the previous subsection shows that 
\[V_0^1, \quad V_4^0 V_0^1 + V_2^0 V_2^1, \quad \cdots, \quad \sum_{l+m=k} V^0_{2l}
V^1_{2m}\]
are spectrally determined. Since each $V^0_k$ is already known (as component of $V_0$), 
the above functions suffice to determine each of $V^1_k$, 
and thus determine $V_1$.  By an inductive 
argument, one can determine all the $V_i$'s. 

To summarize, we have proved:
\begin{theorem}
A semiclassical potential of the form
$\sum_{k \ge 0} \hbar^k V_k\ $ where $V_{0}$ is of the form (\ref{vnot}) with
$a \ne b$, and, for each $k\ V_k$ is analytic and satisfies
$V_{k}(\pm x_1, \pm x_2)=V_k(x_1, x_2)$,
is spectrally determined.
\end{theorem}

\appendix
\section{Proof of theorem \ref{theoremX}}

To prove Theorem \ref{theoremX} 
we must compute the first two non-trivial terms in the full symbol of
the operator $f(S_\chi)W^{l+1}$.

\subsection{The symbol of $f(S_\chi)$.}
Let $\sigma_{S_\chi} \sim H_0 + \h^2 s_2 +\cdots$ be the asymptotic expansion of the symbol
of $S_\chi$.  Note that $s_2 = V + V^{\av}$.
The operator 
\[
U(t) = e^{it S_\chi}
\]
is a pseudodifferential operator of order zero.  Let 
\[
\sigma_U \sim u_0 + \h u_1 + \h^2 u_2+\cdots
\]
be the expansion of its symbol.  Writing down the equation $-i\dot U = S_\chi U$ symbolically 
it is easy to check that $u_0 = e^{itH_0}$, and one computes that
\[
-i\dot{u}_1 = B_1(H_0 , u_0) + H_0 u_1,\quad u_1|_{t=0}= 0,
\]
which has as solution $u_1 = 0$.   $u_2$ is the solution to the problem
\begin{equation}\label{u2dot}
-i\dot{u}_2 = B_2(H_0 , u_0) + H_0 u_2 + u_0 s_2,\quad u_2|_{t=0}= 0.
\end{equation}
Using the general formula
\[
B_2(H, u) = -\frac{1}{4}\sum_{|\alpha|+|\beta|=2} \frac{(-1)^{|\alpha|}}{\alpha!\beta!}\;
\left(\partial^\alpha_x\partial^\beta_p H(x,p)\right)\,\left(\partial^\alpha_p\partial^\beta_x u(x,p)
\right),
\]
one finds:
\begin{equation}\label{}
B_2(H_0 , u_0) = \frac{1}{4}\left(t^2H_0 -int\right)\, e^{itH_0}.
\end{equation}
Substituting into (\ref{u2dot}) and integrating one obtains:
\begin{equation}\label{}
u_2 = ie^{itH_0}\left(\frac{t^3}{12} H_0 + \frac{nt^2}{8i} + t(V+V^{\av})\right).
\end{equation}

We can now prove the following
\begin{lemma}\label{symbolofS}
Let $\sigma_{f(S_\chi)}\sim \phi_0 + \h \phi_1 + \h^2\phi_2+\cdots$ be the asymptotic
expansion of the full symbol of $f(S_\chi)$.  Then one has:
$\phi_0 = f(H_0)$, $\phi_1 = 0$ and
\begin{equation}
\phi_2 = \left( V+ V^{\av}\right)\, f'(H_0) -\frac{n}{8} f''(H_0) - \frac{1}{12}H_0\,f'''(H_0).
\end{equation}
\end{lemma}
\begin{proof}
By the Fourier inversion formula, 
\[
f(S_\chi) =\frac{1}{2\pi}\int U(t)\, \hat{f}(t)\, dt.
\]
At the symbolic level this reads:  $\phi_j = \frac{1}{2\pi}\int u_j(t)\, \hat{f}(t) dt$.
The stated formulae follow from this and the previous calculations of the $u_j$, 
$j=0, 1, 2$.
\end{proof}

\subsection{The end of the proof}

First we need to compute the second non-trivial term in the expansion of the symbol 
of $W^{l+1}$:

\begin{lemma}\label{symbolWpower}
For any non-negative integer $l$ and smooth functions $w_0$, $w_2$ on $\bbR^{2n}$, one 
has:
\[
(w_0 + \h^2 w_2)^{\ast(l+1)} = w_0^{l+1} + 
\h^2\left( (l+1)w_0^l w_2 + \sum_{j=0}^{l-1} w_0^j\, B_2(w_0, w_0^{l-j})\right).
\]
\end{lemma}
The proof is by induction on $l$.

Using lemmas \ref{symbolofS} and \ref{symbolWpower} and recalling that
\[
w_0 = V^{\av}\quad\text{and}\quad w_2 = V^\Delta,
\]
one can easily compute the 
second non-trivial term in the full symbol of $f(S)W^{l+1}$.
The result is precisely the function (\ref{secondinvariant}) plus 
\[
(l+1) f(\frac{|x|^2+|p|^2}{2})(V^{\av})^l V^{\Delta}.
\]
This finishes the proof of Theorem \ref{theoremX}.


\begin{thebibliography}{99}


\bibitem[CV]{CV}
Y. Colin de Verdiere, ``Sur le spectre des op\'erateurs elliptiques \'a bicharat\'eristiques toutes 
\'periodiques", {\sl Comment. Math. Helv.} 54 (1979), 508-522.


\bibitem[Doz]{Doz}
S. Dozias, ``Clustering for the Spectrum of $h$-pseudodifferential Operators with Periodic Flow 
on an Energy Surface", {\sl J. Funct. Anal.} 145 (1997), 296-311.


\bibitem[Gu1]{Gu1}
V. Guillemin, ``Band asymptotics in two dimensions", {\sl Adv. Math. } 42 (1981), 248-282. 

\bibitem[Gu2]{Gu2}
V. Guillemin, ``Spectral Theory on $S^2$: Some Open Questions", {\sl Adv. Math. } 42 (1981), 
283-298. 

\bibitem[GPU]{GPU} V. Gillemin, T. Paul and A. Uribe.  ``'Bottom of the well' semi-classical trace invariants",
{\sl Math. Res. Lett.} 14 (2007), no. 4, 711-719.


\bibitem[GU]{GU}
V. Guillemin and A. Uribe, ``Some Inverse Spectral Results for Semi-classical Schr\"odinger 
Operators", {\sl Math. Res. Lett.} 14 (2007), 623-632.


\bibitem[GUW]{GUW}
V. Guillemin, A. Uribe and Z. Wang, ``A semiclassical heat trace expansion for the perturbed  harmonic oscillator", preprint.  arXiv:1107.2960.

%
\bibitem[GuW]{GuW}
V. Guillemin and Z. Wang, ``Semiclassical Spectral Invariants for Schr\"odinger Operators", {\sl 
submitted.}

%
\bibitem[HR]{HR}
B. Helffer and D. Robert. ``Puits de potentiel g\'en\'eralis\'es et asymptotique semi-
classique", {\sl Ann. Inst. Henri Poincar\'e} 41 (1984), 1-43.


\bibitem[Hez]{Hez}
H. Hezari. ``Inverse Spectral Problems for Schr\"odinger Operators", {\sl Comm. Math. Phys.} 
288 (2009), 1061-1088.


%
\bibitem[Mar]{Mar}
A. Martinez, {\sl An Introduction to Semiclasscial and Microlocal Analysis}, Berlin-Heidelberg-
New York: Springer, 2002. 


\bibitem[O-VV-B]{O-VV-B}  D. Ojeda-Valencia and C. VIllegas Blas,
``On limiting eigenvalue theorems in semiclassical analysis", preprint, 2011.



\bibitem[Ur]{Ur}
A. Uribe, ``Band Invariants and Closed Trajectories on $S^n$", {\sl Adv. Math.} 58 (1985), 
285-299.

\bibitem[UrV]{UrV}
A. Uribe and C. Villegas-Blas, ``Asymptotics of Spectral Clusters for a Perturbation of 
Hydrogen Atom", {\sl Comm. Math. Phys.} 280 (2008), 123-144.

\bibitem[Wei]{Wei}
A. Weinstein, ``Asymptotics of eigenvalue clusters for the Laplacian plus a potential", {\sl Duke 
Math. J.} 44 (1977), 883-892. 

\bibitem[Zel]{Zel}
S. Zelditch, ``Fine Structures of Zoll Spectra", {\sl J. Funct. Anal.} 143 (1997), 415-460.



\end{thebibliography}
\end{document}